\documentclass[10pt,a4paper]{amsart}
\usepackage{graphics,multirow,array,amsmath,amssymb,epsf}

\newtheorem{thm}{Theorem}
\newtheorem{lem}[thm]{Lemma}

\theoremstyle{definition}

\begin{document}
\title{Heptagonal knots and Radon partitions}
\author{Youngsik Huh}
\address{Department of Mathematics,
School of Natural Sciences, Hanyang University, Seoul 133-791,
Korea} \email{yshuh@hanyang.ac.kr}


\thanks{This research was supported by Basic Science Research Program through the National Research Foundation of Korea(NRF) funded by the Ministry of Education, Science and Technology(2010-0009794)}


\begin{abstract}
We establish a necessary and sufficient condition for a heptagonal knot to be figure-8 knot. The condition is described by a set of Radon partitions formed by vertices of the heptagon. In addition we relate this result to the number of nontrivial heptagonal knots in linear embeddings of the complete graph $K_7$ into $\mathbb{R}^3$.
\end{abstract}

\maketitle


\section{Introduction}
An $m$-component {\em link} is a union of $m$ disjoint circles embedded in $\mathbb{R}^3$. Especially a link with only one component is called a {\em knot}. Two knots $K$ and $K^{\prime}$ are said to be {\em ambient isotopic}, denoted by $K\sim K^{\prime}$, if there exists a continuous map $h:\mathbb{R}^3\times[0,1] \rightarrow \mathbb{R}^3$ such that the restriction of $h$ to each $t \in [0,1]$, $h_t:\mathbb{R}^3\times\{t\} \rightarrow \mathbb{R}^3$, is a homeomorphism, $h_0$ is the identity map and $h_1(K_1)=K_2$, to say roughly, $K_1$ can be deformed to $K_2$ without intersecting its strand. The ambient isotopy class of a knot $K$ is called the {\em knot type} of $K$. Especially if $K$ is ambient isotopic to another knot contained in a plane of $\mathbb{R}^3$, then we say that $K$ is {\em trivial}. The ambient isotopy class of links is defined in the same way.

\begin{figure}[h]
\centerline{\epsfxsize=6.5cm \epsfbox{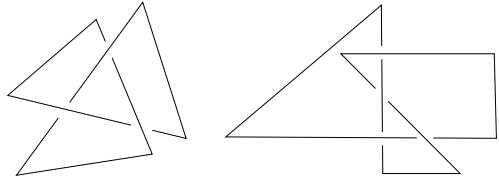}}
\caption{Polygonal presentations of $3_1$ and $4_1$ knots}
\label{fig1}
\end{figure}

In this paper we will focus on polygonal knots. A {\em polygonal knot} is a knot consisting of finitely many line segments, called {\em edges}. The end points of each edge are called {\em vertices}.
Figure \ref{fig1} shows polygonal presentations of two knot types $3_1$ and $4_1$ (These notations for knot types follow the knot tabulation in \cite{Rolfsen}. Usually $3_1$ and its mirror image are called {\em trefoil}, and $4_1$ {\em figure-8}).
For a knot type $\mathfrak{K}$, its {\em polygon index} $p(\mathfrak{K})$ is defined to be the minimal number of edges required to realize $\mathfrak{K}$ as a polygonal knot. Generally it is not easy to determine $p(\mathfrak{K})$ for an arbitrary knot type $\mathfrak{K}$. This quantity was determined only for some specific knot types \cite{Calvo, FLS, J, Mc, Randell-2}. Here we mention a result by Randell on small knots.
\begin{thm} \label{stick-thm}
\cite{Randell-2}
$p(\mbox{trivial knot})=3$, $p(\mbox{trefoil})=6$ and $p(\mbox{figure-8})=7$. Furthermore, $p(\mathfrak{K})\geq 8$ for any other knot type $\mathfrak{K}$.
\end{thm}

Let $V$ be a set of points in $\mathbb{R}^3$. A partition $V_1 \cup V_2$ of $V$ is called a {\em Radon partition} if the two convex hulls of $V_1$ and $V_2$ intersect each other. For example, if $V$ consists of 5 points in general position, then it should have a Radon partition such that $(|V_1|,|V_2|)=(1,4)$ or $(2,3)$.

We remark that the notion of Radon partition can be utilized to describe the knot type of a polygonal knot. In \cite{Alfonsin} a set of Radon partitions are derived from vertices of heptagonal trefoil knots and also hexagonal trefoil knots. Similar work was also done for hexagonal trefoil knot in \cite{Huh}. These results were effectively applied to investigate knots in linear embeddings of the complete graph $K_7$ and $K_6$. An embedding of a graph into $\mathbb{R}^3$ is said to be {\em linear}, if each edge of the graph is mapped to a line segment. In \cite{Alfonsin} Alfons{\'{\i}}n showed that every linear embedding of $K_7$ contains a heptagonal trefoil knot as its cycle. And in \cite{Huh} it was proved that the number of nontrivial knots in any linear embedding of $K_6$ is at most one.

In this paper we give a necessary and sufficient condition for a heptagonal knot to be figure-8 via notion of Radon partition. And we discuss how our result can be utilized to determine the maximal number of heptagonal knots with polygon index 7 residing in linear embeddings of $K_7$.

Now we introduce some notations necessary to describe the main theorem. Let $P$ be a heptagonal knot such that its vertices are in general position. We can label the vertices of $P$ by $\{1, 2, \ldots , 7\}$ so that each vertex $i$ is connected to $i+1$ (mod $7$) by an edge of $P$, that is, a labeling of vertices is determined by a choice of base vertex and an orientation of $P$. Given such a labelling of vertices let $\Delta_{i_1i_2i_3}$ denote the triangle formed by three vertices $\{i_1, i_2, i_3\}$, and $e_{jk}$ the line segment from the vertex $j$ to vertex $k$. The relative position of such a triangle and a line segment will be represented via ``$\epsilon$'' which is defined below:
\begin{itemize}
\item[(i)] If $\Delta_{i_1i_2i_3} \cap e_{jk} = \emptyset$, then set $\epsilon(i_1i_2i_3,jk)=0$.
\item[(ii)] Otherwise,\\
$\epsilon(i_1i_2i_3,jk)= 1$ ({\em resp. $-1$}), when $(\overrightarrow{i_1i_2}\times \overrightarrow{i_2i_3})\cdot \overrightarrow{jk} > 0$ ({\em resp. $<0$}).
\end{itemize}

The tables in Theorem \ref{main-thm} show the values of $\epsilon$ between triangles formed by three consecutive vertices and edges of $P$. If $\epsilon$ is zero, then the corresponding cell in the table is filled by ``$\times$''. Otherwise, we mark by ``$+$'' or ``$-$'' according to the sign of $\epsilon$. For example, according to {\em RS-I},  $\epsilon(123,67)=0$ and \\
\centerline{ $(\epsilon(123,45),\epsilon(123,56),\epsilon(234,56))=(1,-1,-1) \quad \mbox{or} \quad (-1,1,1)$.}
\\
In later sections, for our convenience, we use ``$\bullet$'' to indicate $\epsilon \neq 0$, without specifying the sign.

\begin{thm} \label{main-thm}
Let $P$ be a heptagonal knot such that its vertices are in general position. Then $P$ is figure-8 if and only if the vertices of $P$ can be labelled so that the polygon satisfies one among the three  types RS-I, RS-II and RS-III. \\

\centering{ \small
\begin{tabular}{|c|c|c|c|}
\hline
 & 45 & 56 & 67 \\ \cline{2-4}
\raisebox{1.5ex}[0pt]{123} & $\pm$ & $\mp$ & $\times$ \\ \hline
 & 56 & 67 & 71 \\ \cline{2-4}
\raisebox{1.5ex}[0pt]{234} & $\mp$ & $\times$ & $\times$ \\ \hline
 & 67 & 71 & 12 \\ \cline{2-4}
\raisebox{1.5ex}[0pt]{345} & $\times$ & $\pm$ & $\times$ \\ \hline
 & 71 & 12 & 23 \\ \cline{2-4}
\raisebox{1.5ex}[0pt]{456} & $\pm$ & $\times$ & $\times$ \\ \hline
 & 12 & 23 & 34 \\ \cline{2-4}
\raisebox{1.5ex}[0pt]{567} & $\times$ & $\mp$ & $\times$ \\ \hline
 & 23 & 34 & 45 \\ \cline{2-4}
\raisebox{1.5ex}[0pt]{671} & $\mp$ & $\times$ & $\times$ \\ \hline
 & 34 & 45 & 56 \\ \cline{2-4}
\raisebox{1.5ex}[0pt]{712} & $\times$ & $\pm$ & $\times$ \\ \hline
\multicolumn{4}{c}{ } \\
\multicolumn{4}{c}{RS-I}
\end{tabular}
\hspace{0.5cm}
\begin{tabular}{|c|c|c|c|}
\hline
 & 45 & 56 & 67 \\ \cline{2-4}
\raisebox{1.5ex}[0pt]{123} & $\pm$ & $\mp$ & $\times$ \\ \hline
 & 56 & 67 & 71 \\ \cline{2-4}
\raisebox{1.5ex}[0pt]{234} & $\mp$ & $\times$ & $\times$ \\ \hline
 & 67 & 71 & 12 \\ \cline{2-4}
\raisebox{1.5ex}[0pt]{345} & $\times$ & $\pm$ & $\times$ \\ \hline
 & 71 & 12 & 23 \\ \cline{2-4}
\raisebox{1.5ex}[0pt]{456} & $\pm$ & $\times$ & $\times$ \\ \hline
 & 12 & 23 & 34 \\ \cline{2-4}
\raisebox{1.5ex}[0pt]{567} & $\times$ & $\mp$ & $\times$ \\ \hline
 & 23 & 34 & 45 \\ \cline{2-4}
\raisebox{1.5ex}[0pt]{671} & $\mp$ & $\pm$ & $\times$ \\ \hline
 & 34 & 45 & 56 \\ \cline{2-4}
\raisebox{1.5ex}[0pt]{712} & $\times$ & $\pm$ & $\times$ \\ \hline
\multicolumn{4}{c}{ } \\
\multicolumn{4}{c}{RS-II}
\end{tabular}
\hspace{0.5cm}
\begin{tabular}{|c|c|c|c|}
\hline
 & 45 & 56 & 67 \\ \cline{2-4}
\raisebox{1.5ex}[0pt]{123} & $\pm$ & $\mp$ & $\times$ \\ \hline
 & 56 & 67 & 71 \\ \cline{2-4}
\raisebox{1.5ex}[0pt]{234} & $\times$ & $\mp$ & $\times$ \\ \hline
 & 67 & 71 & 12 \\ \cline{2-4}
\raisebox{1.5ex}[0pt]{345} & $\times$ & $\pm$ & $\times$ \\ \hline
 & 71 & 12 & 23 \\ \cline{2-4}
\raisebox{1.5ex}[0pt]{456} & $\pm$ & $\times$ & $\times$ \\ \hline
 & 12 & 23 & 34 \\ \cline{2-4}
\raisebox{1.5ex}[0pt]{567} & $\times$ & $\mp$ & $\times$ \\ \hline
 & 23 & 34 & 45 \\ \cline{2-4}
\raisebox{1.5ex}[0pt]{671} & $\mp$ & $\times$ & $\times$ \\ \hline
 & 34 & 45 & 56 \\ \cline{2-4}
\raisebox{1.5ex}[0pt]{712} & $\times$ & $\pm$ & $\times$ \\ \hline
\multicolumn{4}{c}{ } \\
\multicolumn{4}{c}{RS-III}
\end{tabular}
}
\end{thm}

In Section 2 we discuss a possible application of Theorem \ref{main-thm}. And the remaining sections will be devoted to the proof of the theorem.

\section{Heptagonal knots in $K_7$ }
In 1983 Conway and Gordon proved that every embedding of $K_7$ into $\mathbb{R}^3$ contains a nontrivial knot as its cycle \cite{CG}. This result was generalized by Negami. He showed that given a knot type $\mathfrak{K}$ there exists a number $r(\mathfrak{K})$ such that every linear embedding of $K_n$ with $n\geq r(\mathfrak{K})$ contains a polygonal knot of type $\mathfrak{K}$ \cite{Negami}.

It would be not easy to determine $r(\mathfrak{K})$ for an arbitrary knot type $\mathfrak{K}$. But if the knot type is of small polygon index, we may attempt to do. For example, Alfons{\'{\i}}n showed that $r(\mbox{trefoil})=7$ \cite{Alfonsin}. To determine the number, he utilized the theory of oriented matroid. This theory provides a way to describe geometric configurations (See \cite{BLSWZ}). Any linear embedding of $K_7$ is determined by fixing the position of seven vertices in $\mathbb{R}^3$. The relative positions of these seven points can be described by an uniform acyclic oriented matroid of rank 4 on seven elements which is in fact a collection of Radon partitions, called {\em signed circuits}, formed by the seven points. Alfons{\'{\i}}n constructed several conditions at least one among which  should be satisfied if a set of seven points constitutes a heptagonal trefoil knot. These conditions are described by a collection of Radon partitions. And then, by help of a computer program, he verified that each of  these matroids satisfies at least one of the conditions. Note that all uniform acyclic oriented matroid of rank 4 on seven elements can be completely listed \cite{Finschi-1, Finschi-2}.

On the other hand we may consider another quantity. Let $\mathcal{F}_n$ be the collection of all linear embeddings of the complete graph $K_n$, and let $c(f)$ be the number of knots with polygon index $n$ in a linear embedding $f \in \mathcal{F}_n$. Define $M(n)$ and $m(n)$ to be
$$ M(n)= \mbox{Max}\; \{c(f) | f \in \mathcal{F}_n \}, \quad m(n)= \mbox{Min}\; \{c(f) | f \in \mathcal{F}_n \}\;. $$
For $n < 6$ these numbers are meaningless because there is no nontrivial knot whose polygonal index is less than $6$. In \cite{Huh} it was shown that $M(6)=1$ and $m(n)=0$ for every $n$. To determine $M(6)$ the author derived a set of Radon partitions from hexagonal trefoil knot. Since $3_1$ and its mirror image are only knot types of polygon index 6, by verifying that the conditional set arises from at most one cycle in any embedded $K_6$, the number $M(6)$ was determined.

We remark that also the number $M(7)$ can be determined by applying our main theorem to a procedure as done by  Alfons{\'{\i}}n. Given a uniform acyclic oriented matroid of rank 4 on seven elements, count the number of permutations which produce any of the conditional partition sets in Theorem \ref{main-thm}. Since the figure-8 is the only knot type of polygon index 7 and the condition in the theorem is necessary and sufficient for a heptagonal knot to be figure-8, the counted number is the number of knots with polygon index 7 in a corresponding embedding of $K_7$. Hence, by getting the maximum among all such numbers over all uniform acyclic oriented matroids of rank 4 on seven elements, $M(7)$ can be determined.

\section{Conway Polynomial}
In this section we give a brief introduction on Conway polynomial which is an ambient isotopy invariant of knots and links. This invariant will be utilized to prove the main theorem in later sections. See \cite{Kunio, Kauffman} for more detailed or kind introduction.

\begin{figure}
\epsfbox{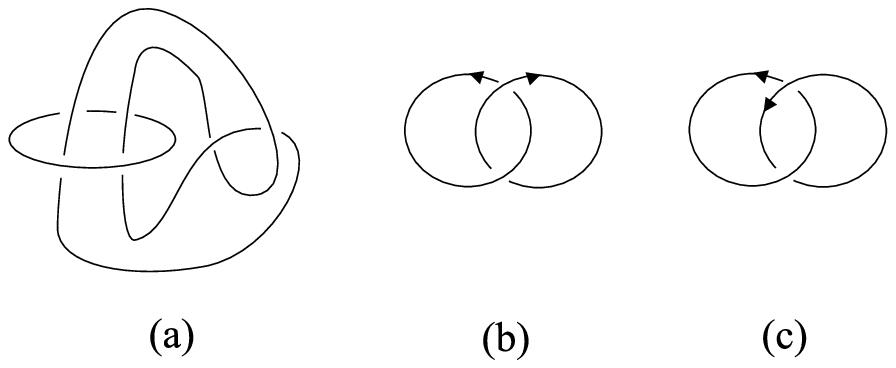}
\caption{}
\label{fig1-1}
\end{figure}

Let $L$ be a link. Given a plane $N$ in $\mathbb{R}^3$, let $\pi_N: N\times \mathbb{R} \rightarrow N$ be the map defined by $\pi(x, y, z)=(x,y)$. Then $\pi_N$ is called a {\em regular projection} of $L$, if the restricted map $\pi_N:L\rightarrow N$ has only finitely many multiple points and every multiple point is a transversal double point. By specifying which strand goes over at each double point of the regular projection, we obtain a {\em diagram} representing $L$. The double points in a diagram are called {\em crossings}. Figure \ref{fig1-1}-(a) shows an example of unoriented link diagram. The diagrams in (b) and (c) represent oriented links. Also the figures in Figure \ref{fig1} can be considered to be unoriented knot diagrams.

Let $L=L_1 \cup L_2$ be a 2-component oriented link. For each $L_i$ we can choose an oriented surface $F_i$ such that $\partial F_i = L_i$. This surface $F_i$ is called a {\em Seifert surface} of $L_i$. Then the linking number $lk(L_1,L_2)$ is defined to be the algebraic intersection number of $L_2$ through $F_1$. It is known that the linking number is independent of the choice of Seifert surface, and $lk(L_1, L_2)=lk(L_2,L_1)$. Hence we may denote the number by  $lk(L)$ instead of $lk(L_1,L_2)$. The linking numbers of the links in Figure \ref{fig1-1}-(b) and (c) are $1$ and $-1$ respectively. The link in (a) is of linking number $0$ for any choice of orientation.

\begin{figure}
\epsfbox{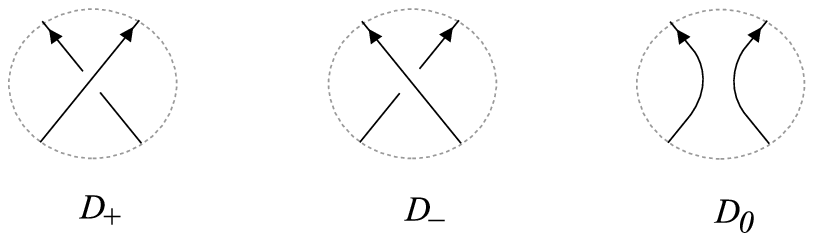}
\caption{}
\label{fig1-2}
\end{figure}

Let $\mathcal{D}$ be the collection of diagrams of all oriented links. Then a function $\nabla : \mathcal{D} \rightarrow \mathbb{Z}[t]$ is uniquely determined by the following three axioms:
\begin{itemize}
\item[(i)] Let $D$ and $D^{'}$ be diagrams which represent two oriented links $L$ and $L^{'}$ respectively. If $L$ is ambient isotopic to $L^{'}$ with orientation preserved, then
$\nabla(D)=\nabla(D^{'})$.
\item[(ii)] If $D$ is a diagram representing the trivial knot, then $\nabla(D)=1$.
\item[(iii)] Let $D_+$, $D_-$ and $D_0$ be three diagrams which are exactly same except at a neighborhood of one crossing point. In the neighborhood they differ as shown in Figure \ref{fig1-2}. The crossing of $D_+$ ({\em resp. $D_-$}) in the figure is said to be {\em positive} ({\em resp. negative}). Then the following equality, called the {\em skein relation}, holds:
$$\nabla(D_+)-\nabla(D_-)=t\nabla(D_0)$$
\end{itemize}
If $D$ is a diagram of an oriented link $L$, then the {\em Conway polynomial} $\nabla(L)$ of $L$ is defined to be $\nabla(D)$. Now we give some facts on Conway polynomial which are necessary for our use in later sections.
\begin{lem} \label{lem-0} \cite{Kunio, Kauffman}
\begin{itemize}
\item[(i)] Let $-K$ be the oriented knot obtained from an oriented knot $K$ by reversing its orientation. Then $\nabla(-K)=\nabla(K)$.
\item[(ii)] If $K$ is trefoil, then $\nabla(K)=1+t^2$. And if $K$ is figure-8, then $\nabla(K)=1-t^2$.
\item[(iii)] Let $L$ be an oriented link with two components. Then its Conway polynomial is of the form $\nabla(L)= a_1t+ a_2t^2+\cdots $ with $a_1 = lk(L)$.
\end{itemize}
\end{lem}
\section{Radon partitions in heptagonal figure-8 knot}
In this section we give several lemmas necessary for the proof of Theorem \ref{main-thm}. Throughout this section $P$ is a heptagonal figure-8 knot such that its vertices are in general position and labelled  by $\{1, 2, \ldots , 7\}$ along an orientation.
Some lemmas will be described by using tables as in Theorem \ref{main-thm}. Note that the blanks in the tables of the following lemma and the rest of this article indicate that the values of $\epsilon$ are not decided yet.

\begin{figure}
\centerline{\epsfbox{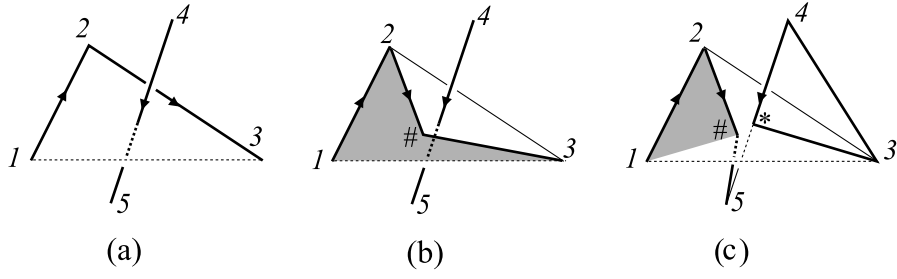}}
\caption{$K_+$, $K_-$ and $K_0$ }
\label{fig2}
\end{figure}

\begin{lem} \label{lem-1}
The following implications hold for $P$.\\

\centering{(i)
{\small
\begin{tabular}{|c|c|c|c|}
\hline
 & 45 & 56 & 67 \\ \cline{2-4}
\raisebox{1.5ex}[0pt]{123} & $\pm$ & $\times$ & $\times$ \\ \hline
\end{tabular}
$\quad \Longrightarrow \quad$
\begin{tabular}{|c|c|c|c|}
\hline
 & 56 & 67 & 71 \\ \cline{2-4}
\raisebox{1.5ex}[0pt]{234} & $\times$ &  &  \\ \hline
 & 67 & 71 & 12 \\ \cline{2-4}
\raisebox{1.5ex}[0pt]{345} & $\mp$ &  &  \\ \hline
\end{tabular} }\\

$\mbox{ }$\\

(ii)
{\small
\begin{tabular}{|c|c|c|c|}
\hline
 & 45 & 56 & 67 \\ \cline{2-4}
\raisebox{1.5ex}[0pt]{123} &$\times$  & $\times$ & $\pm$ \\ \hline
\end{tabular}
$\quad \Longrightarrow \quad$
\begin{tabular}{|c|c|c|c|}
\hline
 & 23 & 34 & 45 \\ \cline{2-4}
\raisebox{1.5ex}[0pt]{671} &  &  & $\mp$ \\ \hline
 & 34 & 45 & 56 \\ \cline{2-4}
\raisebox{1.5ex}[0pt]{712} &  &  & $\times$ \\ \hline
\end{tabular}}
}
\end{lem}
\begin{proof}
Note that (i) is identical with (ii) after relabelling vertices of $P$ along the reverse orientation. Hence it suffices to prove only (i). Assume $\epsilon(123,45)=1$. Then we can choose a diagram of $P$ in which $e_{23}$ and $e_{45}$ produce a positive crossing. Figure \ref{fig2}-(a) depicts the diagram partially. Set $K_+ = P$ and apply the skein relation of Conway polynomial so that
$$\nabla(K_+) - \nabla(K_-) = t\;\nabla(K_0)\; ,$$
where $K_-$ is the cycle $\langle12\#34567\rangle$ and $K_0=\langle12\#567\rangle \cup \langle34*\rangle$ as seen in Figure \ref{fig2}-(b) and (c).
The conditional part of (i) implies that $e_{45}$ is the only edge of $P$ piercing $\Delta_{123}$. Hence $K_- \sim \langle134567\rangle$ by an isotopy in $\Delta_{123}$. Similarly $K_0 \sim \langle1\#567\rangle\cup\langle34*\rangle$.
Since $\langle134567\rangle$ is a hexagon, $K_-$ should be trivial or trefoil by Theorem \ref{stick-thm}. Therefore, $\nabla(K_-)= 1$ or $1+t^2$ and because $\nabla(K_+)= 1-t^2$, we have
$$\nabla(K_0) = -t \quad \mbox{or} \quad -2t \quad .$$
By Lemma \ref{lem-0}-(iii) at least one edge of $\langle1\#567\rangle$ penetrates $\Delta_{34*}$ in negative direction. Note that $\Delta_{34*}$ is contained in a half space $H_{123}^-$ with respect to the plane $H_{123}^0$ formed by $\{1, 2, 3 \}$. Since $e_{1\#}$ belongs to $\Delta_{123}$ and $e_{\#5}$ belongs to another half space $H_{123}^+$, the two edges are excluded from candidates. Also $e_{56}$ is excluded because $\Delta_{34*}\subset\Delta_{345}$. Hence $e_{67}$ and $e_{71}$ are the only edges which may penetrate $\Delta_{34*}$.
But the vertex $1$ belongs to $H_{34*}^+$, which implies that if $e_{71}$ penetrates $\Delta_{34*}$, then the orientation of intersection should be positive. Therefore we can conclude $\epsilon(34*,67)=-1$, and hence $\epsilon(345,67)=-1$.

Let $T_{5,123}^{\infty}$ be the set of all half infinite lines starting the vertex $5$ and passing through a point of $\Delta_{123}$. Clearly $\Delta_{234} \subset T_{5,123}^{\infty}$. Hence if we suppose $\epsilon(234,56)\neq 0$, then also $\epsilon(123,56)\neq 0$, which is contradictory to the condition of (i).

In the case that $\epsilon(123,45)=-1$ we can prove the implication in a similar way.
\end{proof}

\begin{figure}
\centerline{\epsfbox{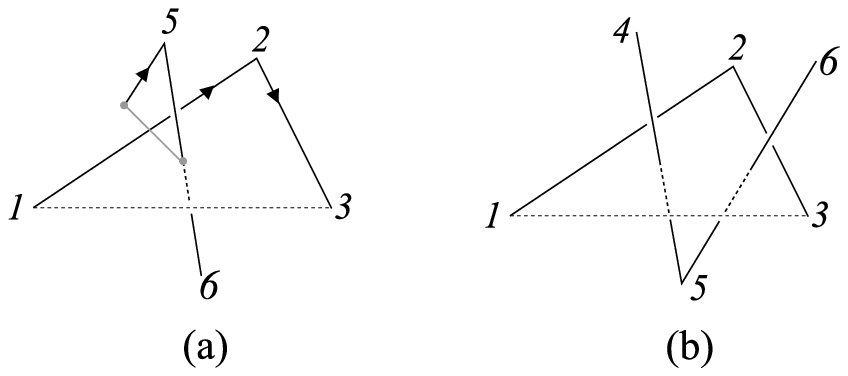}}
\caption{ }
\label{fig3}
\end{figure}

\begin{lem} \label{lem-2}
$$\mbox{(i)} \quad \epsilon(123,56)=\pm1 \;\; \mbox{and} \;\; \epsilon(456,12)\neq 0 \quad \Longrightarrow \quad \epsilon(456,12)=\pm1 $$
$$\mbox{(ii)} \quad \epsilon(123,56)=\pm1 \;\; \mbox{and} \;\; \epsilon(567,23)\neq 0 \quad \Longrightarrow \quad  \epsilon(567,23)=\pm1$$
\end{lem}
\begin{proof}
Assuming $\epsilon(123,56)=1$, the conditional part of (i) can be illustrated as Figure \ref{fig3}-(a). From the figure we can verify (i). Similarly (ii) can be proved.
\end{proof}

\begin{lem} \label{lem-3} $\mbox{ }$\\

\centering{(i)
{\small
\begin{tabular}{|c|c|c|c|}
\hline
 & 45 & 56 & 67 \\ \cline{2-4}
\raisebox{1.5ex}[0pt]{123} & $\bullet$ & $\bullet$ &  \\ \hline
\end{tabular}
$\quad \Longrightarrow \quad$
\begin{tabular}{|c|c|c|c|}
\hline
 & 67 & 71 & 12 \\ \cline{2-4}
\raisebox{1.5ex}[0pt]{345} &  &  & $\times$ \\ \hline
 & 71 & 12 & 23 \\ \cline{2-4}
\raisebox{1.5ex}[0pt]{456} & $\bullet$ & $\times$ & $\times$ \\ \hline
\end{tabular} }\\

$\mbox{ }$\\

(ii)
{\small
\begin{tabular}{|c|c|c|c|}
\hline
 & 45 & 56 & 67 \\ \cline{2-4}
\raisebox{1.5ex}[0pt]{123} &  & $\bullet$ & $\bullet$ \\ \hline
\end{tabular}
$\quad \Longrightarrow \quad$
\begin{tabular}{|c|c|c|c|}
\hline
 & 12 & 23 & 34 \\ \cline{2-4}
\raisebox{1.5ex}[0pt]{567} & $\times$ & $\times$ & $\bullet$ \\ \hline
 & 23 & 34 & 45 \\ \cline{2-4}
\raisebox{1.5ex}[0pt]{671} & $\times$ &  &  \\ \hline
\end{tabular}}}
\end{lem}
\begin{proof}
Assuming $\epsilon(123,45)=1$, the conditional part of (i) can be illustrated as Figure \ref{fig3}-(b). The figure clearly shows that $\epsilon(345,12)=0$, $\epsilon(456,12)=0$ and $\epsilon(456,23)=0$.
Note that $e_{71}$, $e_{12}$ and $e_{23}$ are the only possible edges of $P$ which may penetrate $\Delta_{456}$. Hence $\epsilon(456,71)$ should be nonzero. Otherwise, $P=\langle1234567\rangle$ can be isotoped to the hexagon $\langle123467\rangle$ along $\Delta_{456}$, which contradicts that $P$ is of polygon index 7 by Theorem \ref{stick-thm}.
Similarly (ii) can be proved.
\end{proof}

\begin{figure}
\centerline{\epsfbox{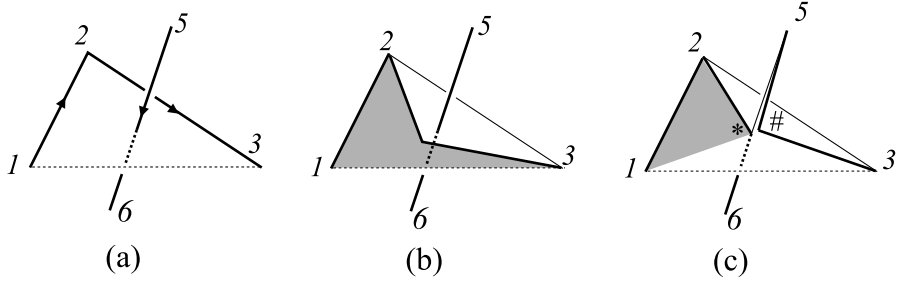}}
\caption{$K_+$, $K_-$ and $K_0$ }
\label{fig4}
\end{figure}

\begin{lem} \label{lem-4}
$P$ does not allow any of two cases below:\\

\centering{
(i) {\small
\begin{tabular}{|c|c|c|c|}
\hline
 & 45 & 56 & 67 \\ \cline{2-4}
\raisebox{1.5ex}[0pt]{123} & $\times$ & $\pm$ & $\times$ \\ \hline
 & 67 & 71 & 12 \\ \cline{2-4}
\raisebox{1.5ex}[0pt]{345} & $\times$ & $\pm$ &  \\ \hline
\end{tabular} }
(ii) {\small
\begin{tabular}{|c|c|c|c|}
\hline
 & 23 & 34 & 45 \\ \cline{2-4}
\raisebox{1.5ex}[0pt]{671} &  & $\pm$ & $\times$ \\ \hline
 & 45 & 56 & 67 \\ \cline{2-4}
\raisebox{1.5ex}[0pt]{123} & $\times$ & $\pm$ & $\times$ \\ \hline
\end{tabular} }
}
\end{lem}
\begin{proof}
Suppose that (i) is true. It suffices to consider the case that $\epsilon(123,56)=\epsilon(345,71)=1$.
Apply the skein relation to the crossing between $e_{23}$ and $e_{56}$ as seen in Figure \ref{fig4}, so that $K_- \sim \langle134567\rangle$ and $K_0 \sim \langle1*67\rangle \cup \langle5\#34\rangle$. Then $\nabla(K_0)$ should be $-t$ or $-2t$, that is, the linking number of $K_0$ is $-1$ or $-2$. Note that $\Delta_{5\#3}\cup\Delta_{345}$ is a Seifert surface of $\langle5\#34\rangle$. Therefore
$$\sum_{i \in \{1*, *6, 67, 71\} }(\epsilon(5
\#3,e_i) + \epsilon(345,e_i))= -1 \;\; \mbox{or} \;\; -2 \;\;.$$
By our assumption $\epsilon(345,67)=0$ and $\epsilon(345,71)=1$. Clearly we know that $\epsilon(5\#3,e_i)=0$ for $i=1*, *6$.
Also $\epsilon(345,*6)$ is $0$ because $e_{*6}$ is a segment of $e_{56}$. Select the point $\#$ so that $\Delta_{5\#3} \subset \Delta_{356}$, hence  $\epsilon(5\#3, 67)$ is $0$.
Therefore the summation should be $-1$, and $\epsilon(345,1*)=\epsilon(5\#3,71)=-1$. But the vertex $1$ belongs to $H_{5\#3}^{+}$ as seen in the figure, hence if $e_{71}$ penetrates $\Delta_{5\#3}$, then the orientation of intersection should be positive, which is a contradiction.

(ii) is derived directly from (i) by relabelling the vertices after reversing the orientation of $P$.
\end{proof}

\begin{figure}
\centerline{\epsfbox{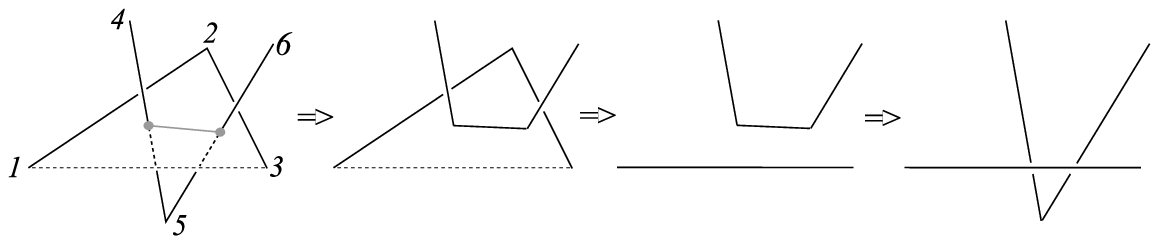}}
\caption{ }
\label{fig4-1}
\end{figure}

Two integers $i$ and $j$ indicate the same vertex if $i \equiv j$ (mod 7). For an integer $i$, define $I(i)$ to be the number of edges of $P$ penetrating $\Delta_{i,i+1,i+2}$, that is,
$$I(i) = \sum_{j \in \{ i+3, i+4, i+5\}} |\epsilon((i,i+1,i+2),( j,j+1))| \;\;.$$
\begin{lem} \label{lem-5}
There exists no integer $i$ such that $I(i) \geq 2$ and $I(i+1) \geq 2$.
\end{lem}
\begin{proof} Suppose that $I(1) \geq 2$ and $I(2)\geq2$. Then, up to the relabelling \\
\centerline{$(1,2,3,4,5,6,7)\rightarrow(4,3,2,1,7,6,5)$}\\
it is enough to observe the following six cases:

\centerline{
(i) {\small
\begin{tabular}{|c|c|c|c|}
\hline
 & 45 & 56 & 67 \\ \cline{2-4}
\raisebox{1.5ex}[0pt]{123} & $\bullet$ & $\bullet$ & \\ \hline
 & 56 & 67 & 71 \\ \cline{2-4}
\raisebox{1.5ex}[0pt]{234} & $\bullet$ & $\bullet$ & \\ \hline
\end{tabular} } \hspace{0.5cm}
(ii) {\small
\begin{tabular}{|c|c|c|c|}
\hline
 & 45 & 56 & 67 \\ \cline{2-4}
\raisebox{1.5ex}[0pt]{123} & $\bullet$ & $\bullet$ & \\ \hline
 & 56 & 67 & 71 \\ \cline{2-4}
\raisebox{1.5ex}[0pt]{234} & $\bullet$ &  & $\bullet$\\ \hline
\end{tabular} } \hspace{0.5cm}
(iii) {\small
\begin{tabular}{|c|c|c|c|}
\hline
 & 45 & 56 & 67 \\ \cline{2-4}
\raisebox{1.5ex}[0pt]{123} & $\bullet$ & $\bullet$ & \\ \hline
 & 56 & 67 & 71 \\ \cline{2-4}
\raisebox{1.5ex}[0pt]{234} &  & $\bullet$ & $\bullet$\\ \hline
\end{tabular}}
}
\centerline{
(iv) {\small
\begin{tabular}{|c|c|c|c|}
\hline
 & 45 & 56 & 67 \\ \cline{2-4}
\raisebox{1.5ex}[0pt]{123} & $\bullet$ &  &$\bullet$ \\ \hline
 & 56 & 67 & 71 \\ \cline{2-4}
\raisebox{1.5ex}[0pt]{234} & $\bullet$ &$\bullet$  & \\ \hline
\end{tabular} } \hspace{0.5cm}
(v) {\small
\begin{tabular}{|c|c|c|c|}
\hline
 & 45 & 56 & 67 \\ \cline{2-4}
\raisebox{1.5ex}[0pt]{123} & $\bullet$ &  &$\bullet$ \\ \hline
 & 56 & 67 & 71 \\ \cline{2-4}
\raisebox{1.5ex}[0pt]{234} & $\bullet$ &  & $\bullet$\\ \hline
\end{tabular} } \hspace{0.5cm}
(vi) {\small
\begin{tabular}{|c|c|c|c|}
\hline
 & 45 & 56 & 67 \\ \cline{2-4}
\raisebox{1.5ex}[0pt]{123} &  & $\bullet$ &$\bullet$ \\ \hline
 & 56 & 67 & 71 \\ \cline{2-4}
\raisebox{1.5ex}[0pt]{234} & $\bullet$ &$\bullet$  & \\ \hline
\end{tabular} }
}

\vspace{0.1cm}
For Case (i), apply Lemma \ref{lem-3} to $\Delta_{123}$ and $\Delta_{234}$. Then we have

\centerline{
\small
\begin{tabular}{|c|c|c|c|}
\hline
 & 71 & 12 & 23 \\ \cline{2-4}
\raisebox{1.5ex}[0pt]{456} & $\bullet$ & $\times$ & $\times$ \\ \hline
 & 12 & 23 & 34 \\ \cline{2-4}
\raisebox{1.5ex}[0pt]{567} & $\bullet$ & $\times$ & $\times$ \\ \hline
\end{tabular}
}
\noindent But, applying Lemma \ref{lem-1}-(i) to $\Delta_{456}$, $\epsilon(567,12)$ should be $0$, a contradiction. Also Case (iii) can be rejected in a similar way. For Case (vi) to be excluded, only Lemma \ref{lem-3} is enough.

For Case (ii) we may assume further that $\epsilon(123,45)=1$. Then, for $e_{71}$ to penetrate $\Delta_{234}$, the vertex $7$ should belong to $H_{123}^-$. This implies that the region $\Delta_{456}\cap H_{123}^+$ is not penetrated by any edge of $P$, and $e_{45}$ and $e_{56}$ are the only edges of $P$ penetrating $\Delta_{123}$. Hence we can isotope $P$ to $\langle134567\rangle$ as illustrated in Figure \ref{fig4-1}, which contradicts Theorem \ref{stick-thm}. Also in Case (iv) $P$ can be isotoped to $\langle245671\rangle$ in a similar way.

For case (v), we may suppose further that $\epsilon(123,45)=1$. Then $\Delta_{234}$ belongs to $H_{123}^-$. If $\epsilon(123,67)=1$, then $e_{71}$ belongs to $H_{123}^+$ and hence can not penetrate $\Delta_{234}$, a contradiction. Similarly if $\epsilon(123,67)=-1$, then $e_{56}$ can not penetrate $\Delta_{234}$.
\end{proof}

\begin{figure}
\centerline{\epsfbox{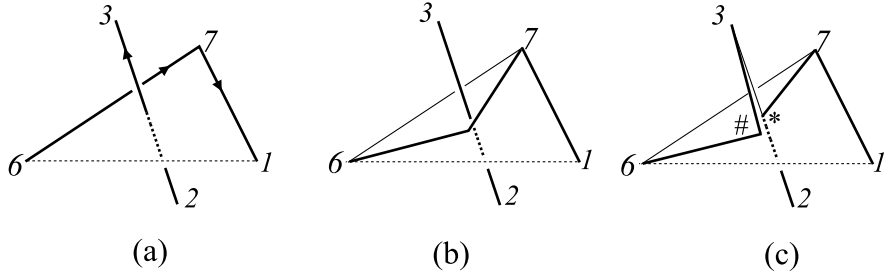}}
\caption{$K_-$, $K_+$ and $K_0$ }
\label{fig5}
\end{figure}

\begin{lem} \label{lem-6}
There exists no pair of distinct integers $(i,j)$ such that

\centerline{
{\small
\begin{tabular}{|c|c|c|c|}
\hline
 & i+3, i+4 & i+4, i+5 & i+5, i+6 \\ \cline{2-4}
\raisebox{1.5ex}[0pt]{i, i+1, i+2} & $\bullet$ & $\bullet$ & $\times$\\ \hline
 & j+3, j+4 & j+4, j+5 & j+5, j+6 \\ \cline{2-4}
\raisebox{1.5ex}[0pt]{j, j+1, j+2} & $\times$ & $\bullet$ & $\bullet$\\ \hline
\end{tabular} }
}
\end{lem}
\begin{proof} By Lemma \ref{lem-5} it is enough to observe the four cases: $(i,j)=(1,3)$, $(1,4)$, $(1,5)$ and $(1,6)$.
The first two cases are contradictory to Lemma \ref{lem-3}. For the fourth case, applying Lemma \ref{lem-3} to $\Delta_{123}$, we have

\centerline{\small
\begin{tabular}{|c|c|c|c|}
\hline
 & 71 & 12 & 23 \\ \cline{2-4}
\raisebox{1.5ex}[0pt]{456} & $\bullet$ & $\times$ & $\times$\\ \hline
\end{tabular} }

\noindent And apply Lemma \ref{lem-1} to $\Delta_{456}$, to have $\epsilon(671,23)\neq 0$, a contradiction.

Lastly suppose $(i,j)=(1,5)$. Then, we can observe which edges penetrate $\Delta_{671}$ and $\Delta_{712}$ as follows:

{\small
\begin{tabular}{|c|c|c|c|}
\hline
 & 45 & 56 & 67 \\ \cline{2-4}
\raisebox{1.5ex}[0pt]{123} & $\bullet$ & $\bullet$ & $\times$\\
\hline
 & 12 & 23 & 34 \\ \cline{2-4}
\raisebox{1.5ex}[0pt]{567} & $\times$ & $\bullet$ & $\bullet$\\ \hline
\end{tabular}
\begin{tabular}{c}
Lemma \ref{lem-3}\\ $\Longrightarrow$
\end{tabular}
\begin{tabular}{|c|c|c|c|}
\hline
 & 56 & 67 & 71 \\ \cline{2-4}
\raisebox{1.5ex}[0pt]{234} & $\times$ & $\times$ & $\bullet$\\
\hline
 & 71 & 12 & 23 \\ \cline{2-4}
\raisebox{1.5ex}[0pt]{456} & $\bullet$ & $\times$ & $\times$\\ \hline
\end{tabular}
\begin{tabular}{c}
Lemma \ref{lem-1}\\ $\Longrightarrow$
\end{tabular}
\begin{tabular}{|c|c|c|c|}
\hline
 & 23 & 34 & 45 \\ \cline{2-4}
\raisebox{1.5ex}[0pt]{671} & $\bullet$ &  & \\
\hline
 & 34 & 45 & 56 \\ \cline{2-4}
\raisebox{1.5ex}[0pt]{712} &  &  & $\bullet$\\ \hline
\end{tabular}
\begin{tabular}{c}
Lemma \ref{lem-5}\\ $\Longrightarrow$ \\ $(i,j)=(1,5)$
\end{tabular}
\begin{tabular}{|c|c|c|c|}
\hline
 & 23 & 34 & 45 \\ \cline{2-4}
\raisebox{1.5ex}[0pt]{671} & $\bullet$ & $\times$ & $\times$\\
\hline
 & 34 & 45 & 56 \\ \cline{2-4}
\raisebox{1.5ex}[0pt]{712} & $\times$ & $\times$ & $\bullet$\\ \hline
\end{tabular}
}

\noindent We may assume $\epsilon(123,45)=1$. Then clearly $\epsilon(123,56)=-1$, and  by Lemma \ref{lem-1} $\epsilon(671,23)=-1$. Now we apply the skein relation to $e_{23} \cup e_{67}$ as seen in Figure \ref{fig5} so that $K_- = P$, $K_0 \sim \langle6\#345\rangle \cup \langle2*1\rangle$ and $\Delta_{2*1} \subset \Delta_{123}$. Then immediately it is observed that $$\epsilon(2*1,6\#)=\epsilon(2*1,\#3)=\epsilon(2*1,34)=0$$
Recall $\epsilon(123,56)=-1$, which implies that if $\epsilon(2*1,56) \neq 0$, then the value should be negative. Therefore, considering $\nabla(K_0)=t$ or $2t$, it should hold that
$$\epsilon(2*1,45)= 1, \quad \epsilon(2*1,56)= 0$$
and $e_{56}$ penetrates $\Delta_{*31}=\Delta_{123}-\Delta_{2*1}$.
Since $\Delta_{*31} \subset H_{671}^-$, the vertex $5$ belongs to $H_{671}^-$. Hence $e_{56}\subset H_{671}^-$, which is contradictory to $\epsilon(712,56)\neq0$ because $\Delta_{712}$ belongs to the other half space $H_{671}^+$.
\end{proof}

\begin{lem} \label{lem-7}
$\mbox{ }$

(i) For every $i$, $I(i) \geq 1$.

(ii) There exists an integer $i$ such that $I(i) \geq 2$.

(iii) For every $i$, $I(i) <3$.

(iv) If $I(i)=2$ for some $i$, then $e_{i+4,i+5}$ should penetrate $\Delta_{i,i+1,i+2}$.
\end{lem}
\begin{proof}[Proof of (i).] Suppose $I(1)=0$, that is, $\Delta_{123}$ is not penetrated by any edge of $P$. Then we can isotope $P$ along $\Delta_{123}$, so that $P \sim \langle134567\rangle$. By Theorem \ref{stick-thm}, the hexagon $\langle134567\rangle$ is trivial or trefoil, a contradiction.
\end{proof}

\begin{proof}[Proof of (ii).] Suppose  $I(i)=1$ for every $i$. Then, among $e_{45}$, $e_{56}$ and $e_{67}$, only one edge penetrates $\Delta_{123}$. Firstly assume that $e_{45}$ does. Then, applying Lemma \ref{lem-1} repeatedly, we have a sequence of implications:

\centerline{\small
\begin{tabular}{|c|c|c|c|}
\hline
 & 45 & 56 & 67 \\ \cline{2-4}
\raisebox{1.5ex}[0pt]{123} & $\bullet$ &  $\times$ & $\times$\\
\hline
\end{tabular}
$\Longrightarrow$
\begin{tabular}{|c|c|c|c|}
\hline
 & 67 & 71 & 12 \\ \cline{2-4}
\raisebox{1.5ex}[0pt]{345} & $\bullet$ &  $\times$ & $\times$\\
\hline
\end{tabular}
$\Longrightarrow$
\begin{tabular}{|c|c|c|c|}
\hline
 & 12 & 23 & 34 \\ \cline{2-4}
\raisebox{1.5ex}[0pt]{567} & $\bullet$ &  $\times$ & $\times$\\
\hline
\end{tabular}
$\Longrightarrow$
}

\centerline{\small
\begin{tabular}{|c|c|c|c|}
\hline
 & 34 & 45 & 56 \\ \cline{2-4}
\raisebox{1.5ex}[0pt]{712} & $\bullet$ &  $\times$ & $\times$\\
\hline
\end{tabular}
$\Longrightarrow$
\begin{tabular}{|c|c|c|c|}
\hline
 & 56 & 67 & 71 \\ \cline{2-4}
\raisebox{1.5ex}[0pt]{234} & $\bullet$ &  $\times$ & $\times$\\
\hline
\end{tabular}
}

\noindent But, by Lemma \ref{lem-1} again, the first and last tables are contradictory to each other. The case that $e_{67}$ penetrates the triangle is rejected in a similar way.

Now it can be assumed that every $\Delta_{i,i+1,i+2}$ is penetrated only by $e_{i+4,i+5}$. Then, applying Lemma \ref{lem-2} repeatedly, we have that

\centerline{\small
\begin{tabular}{|c|c|c|c|}
\hline
 & 45 & 56 & 67 \\ \cline{2-4}
\raisebox{1.5ex}[0pt]{123} &  $\times$ &  $\pm$ & $\times$\\
\hline
\end{tabular}
 and
\begin{tabular}{|c|c|c|c|}
\hline
 & 67 & 71 & 12 \\ \cline{2-4}
\raisebox{1.5ex}[0pt]{345} &  $\times$ &  $\pm$ & $\times$\\
\hline
\end{tabular} ,
}

\noindent which is contradictory to Lemma \ref{lem-4}.
\end{proof}
\begin{proof}[Proof of (iii).]
Suppose $I(1)=3$. Then we have two implications as follows:
\begin{flushleft}\small
\begin{tabular}{|c|c|c|c|}
\hline
 & 45 & 56 & 67 \\ \cline{2-4}
\raisebox{1.5ex}[0pt]{123} &  $\bullet$ &  $\bullet$ & \\
\hline
\end{tabular}
\begin{tabular}{c}
Lemma \ref{lem-3}\\ $\Longrightarrow$
\end{tabular}
\begin{tabular}{|c|c|c|c|}
\hline
 & 71 & 12 & 23 \\ \cline{2-4}
\raisebox{1.5ex}[0pt]{456} &  $\bullet$ &  $\times$ & $\times$ \\
\hline
\end{tabular}
\begin{tabular}{c}
Lemma \ref{lem-1}\\ $\Longrightarrow$
\end{tabular}
\begin{tabular}{|c|c|c|c|}
\hline
 & 23 & 34 & 45 \\ \cline{2-4}
\raisebox{1.5ex}[0pt]{671} &  $\bullet$ &  & \\
\hline
\end{tabular}
\end{flushleft}
\begin{flushleft} \small
\begin{tabular}{|c|c|c|c|}
\hline
 & 45 & 56 & 67 \\ \cline{2-4}
\raisebox{1.5ex}[0pt]{123} &   &  $\bullet$ & $\bullet$ \\
\hline
\end{tabular}
\begin{tabular}{c}
Lemma \ref{lem-3}\\ $\Longrightarrow$
\end{tabular}
\begin{tabular}{|c|c|c|c|}
\hline
 & 23 & 34 & 45 \\ \cline{2-4}
\raisebox{1.5ex}[0pt]{671} &  $\times$ &  &  \\
\hline
\end{tabular}
\end{flushleft}

\noindent But these are contradictory to each other.
\end{proof}

\begin{figure}
\centerline{\epsfbox{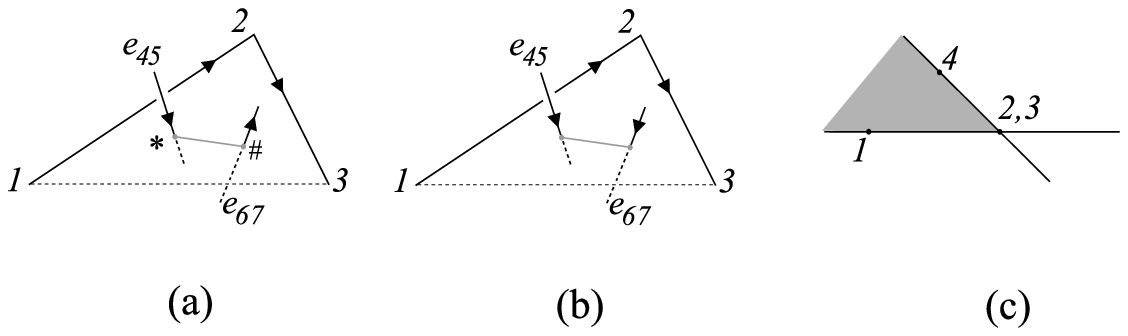}}
\caption{ }
\label{fig7}
\end{figure}

\begin{proof}[Proof of (iv).]
Suppose that $\Delta_{123}$ satisfies the following:

\centerline{\small
\begin{tabular}{|c|c|c|c|}
\hline
 & 45 & 56 & 67 \\ \cline{2-4}
\raisebox{1.5ex}[0pt]{123} & $\bullet$  &  $\times$ & $\bullet$ \\
\hline
\end{tabular}}

\noindent Then it is enough to observe two cases $(\epsilon(123,45), \epsilon(123,67))=(1, -1)$ and $(1, 1)$. These cases are depicted as  in Figure \ref{fig7}-(a) and (b) respectively. In the first case $5$ and $6$ are the only vertices which belong to $H_{123}^+$. Therefore P can be isotoped to $\langle1234*\#7\rangle$ along the tetragon formed by $\{ *, 5, 6, \#\}$. And lift $e_{*\#}$ slightly into $H_{123}^-$, then we have $I(1)=0$ for the resulting heptagon, a contradiction.

For the second case we observe which edges penetrate $\Delta_{234}$. Note that $\Delta_{234} \subset T_{5,123}^{\infty}$. Hence if a line starting at the vertex $5$ penetrates $\Delta_{234}$, then it also penetrates $\Delta_{123}$. This implies $\epsilon(234,56)=0$. Also $\epsilon(234,71)=0$, because $e_{71}$ belongs to $H_{123}^+$ but $\Delta_{234}$ belongs to the other half space $H_{123}^-$. Therefore $e_{67}$ is the only edge of $P$ penetrating $\Delta_{234}$. Furthermore the orientation of intersection should be positive. This can be seen easily from Figure \ref{fig7}-(c). Let $N$ be a plane in $\mathbb{R}^3$ orthogonal to $\overrightarrow{23}$. And let $\pi : \mathbb{R}^3 \equiv N\times\mathbb{R} \rightarrow N$ be the orthogonal projection onto $N$ such that the vertex $3$ is above the vertex $2$ with respect to the $\mathbb{R}$-coordinate. Figure \ref{fig7}-(c) depicts the image of $H_{123}^0 \cup H_{234}^0$ under $\pi$. Suppose $\epsilon(234,67)=-1$. Since $\epsilon(123,67)=1$, the vertex 6 should belong to $H_{123}^- \cap H_{234}^+$ which corresponds to the shaded region in the figure. Then, as seen in the figure, it is impossible that $e_{67}$ penetrates both  $\Delta_{123}$ and $\Delta_{234}$.

In a similar way $e_{45}$ should be the only edge of $P$ penetrating $\Delta_{712}$ and the orientation of intersection is positive. To summarize, we have

\centerline{ {\small
\begin{tabular}{|c|c|c|c|}
\hline
 & 56 & 67 & 71 \\ \cline{2-4}
\raisebox{1.5ex}[0pt]{234} & $\times$  &  $+$ & $\times$ \\
\hline
\end{tabular}}
 and
{\small
\begin{tabular}{|c|c|c|c|}
\hline
 & 34 & 45 & 56 \\ \cline{2-4}
\raisebox{1.5ex}[0pt]{712} & $\times$  &  $+$ & $\times$ \\
\hline
\end{tabular}}}

\noindent This contradicts Lemma \ref{lem-4}.
\end{proof}

\section{Proof of theorem \ref{main-thm}}
We prove the ``only if'' part of Theorem \ref{main-thm} by filling in the table of penetrations in $P$. By Lemma \ref{lem-7} it can be assumed that

\centerline{ \small
\begin{tabular}{|c|c|c|c|}
\hline
 & 45 & 56 & 67 \\ \cline{2-4}
\raisebox{1.5ex}[0pt]{123} & $\bullet$  &  $\bullet$ & $\times$ \\
\hline
\end{tabular}}

Applying Lemma \ref{lem-3} to $\Delta_{123}$ and Lemma \ref{lem-1} to $\Delta_{456}$, we have the initial status $S_0$ as shown in Table \ref{table-1}. Considering Lemma \ref{lem-6}, we know that the row of $567$ should be filled by $(\times,\bullet,\times)$ or $(\times,\times,\bullet)$. The second is excluded by Lemma \ref{lem-1}. Lemma \ref{lem-7}-(iv) guarantees $\epsilon(671,45)=0$. Hence the status $S^{\prime}_{0}$ is derived. Observe how the row of $712$ can be filled. $I(7)$ should be $1$ or $2$ by Lemma \ref{lem-7}-(i) and (iii). In fact $I(7)$ should be 1 by Lemma \ref{lem-5}. And $(\bullet,\times,\times)$ and $(\times,\times,\bullet)$ are disallowed by Lemma \ref{lem-1}, hence $S^{\prime\prime}_{0}$ is derived.

In a similar way we know that $234$ should have $(\bullet,\times,\times)$ or $(\times,\bullet, \times)$. Hence $S^{\prime\prime}_{0}$ can proceed to the status $S_1$ or $S_2$.

\noindent \emph{Case 1: the row of $234$ is filled with $(\bullet,\times,\times)$.} Observe $345$ in $S_1$. Apply Lemma \ref{lem-1} to $\Delta_{234}$. Then $\epsilon(345,67)=0$, and $\epsilon(345,71)\neq 0$. Therefore we obtain $S_{1}^{\prime}$.
Finally, let $S_{1-1}^{\prime}$ ({\em resp.} $S_{1-2}^{\prime}$) be the status obtained from $S_{1}^{\prime}$ by setting $\epsilon(671,34)$ to be zero ({{\em resp.} nonzero}).

Note that if the table is completely filled with ``$\bullet$'' and``$\times$'', then the orientation of intersection is automatically determined. See $S_{1-1}^{\prime}$. Since $123$ has $(\bullet,\bullet,\times)$, the possible orientation is $(+,-,\times)$ or $(-,+,\times)$. Assume the former. Applying Lemma \ref{lem-1} to $\Delta_{671}$, we know $\epsilon(671,23)=-1$. Also applying the lemma to $\Delta_{456}$ and $\Delta_{234}$, we have $\epsilon(456,71)=1$ and $\epsilon(234,56)=-1$. Furthermore, from the assumption $\epsilon(123,56)=-1$, it is derived that $\epsilon(567,23)=-1$ and $\epsilon(345,71)=1$   by Lemmas \ref{lem-2} and \ref{lem-4} respectively. Similarly $\epsilon(712,45)=1$. Therefore $S_{1-1}^{\prime}$ is identical with {\em RS-I}.

For $S_{1-2}^{\prime}$, first determine the orientations in the second column following the method used above. Then, under the assumption $\epsilon(123,56)=-1$, it should hold that $\epsilon(671,34)=1$. This implies $\epsilon(671,23)=-1$, from which the orientations in the first column can be determined. In this way we can verify that $S_{1-2}^{\prime}$ is identical with {\em RS-II}.

\emph{Case 2: the row of $234$ is filled with $(\times,\bullet,\times)$.} If $671$ has $(\bullet,\bullet,\times)$ in $S_2$, then $234$ should have $(\bullet,\times,\times)$ by Lemma \ref{lem-3}, a contradiction. Therefore $S_2$ proceeds only to $S_2^{\prime}$. Suppose that $345$ can be filled with $(\bullet,\times,\times)$. Then also we can derive a contradiction by applying Lemma \ref{lem-1} to $\Delta_{345}$. Hence, in $S_2^{\prime}$, $345$ should have $(\bullet,\bullet,\times)$ or $(\times,\bullet,\times)$. In the former case we have $S_{2-1}^{\prime}$ which becomes {\em RS-II} after relabelling vertices by the cyclic permutation sending $(3,4,5)$ to $(1,2,3)$. In the latter we have $S_{2-2}^{\prime}$ which is {\em RS-III}.
\begin{table}
{ \tiny
\begin{tabular}{|c|c|c|c|}
\hline
 & 45 & 56 & 67 \\ \cline{2-4}
\raisebox{1.5ex}[0pt]{123} & $\bullet$ & $\bullet$ & $\times$ \\ \hline
 & 56 & 67 & 71 \\ \cline{2-4}
\raisebox{1.5ex}[0pt]{234} & &  &  \\ \hline
 & 67 & 71 & 12 \\ \cline{2-4}
\raisebox{1.5ex}[0pt]{345} & &  & $\times$ \\ \hline
 & 71 & 12 & 23 \\ \cline{2-4}
\raisebox{1.5ex}[0pt]{456} & $\bullet$ & $\times$ & $\times$ \\ \hline
 & 12 & 23 & 34 \\ \cline{2-4}
\raisebox{1.5ex}[0pt]{567} & $\times$ &  & \\ \hline
 & 23 & 34 & 45 \\ \cline{2-4}
\raisebox{1.5ex}[0pt]{671} & $\bullet$ & &  \\ \hline
 & 34 & 45 & 56 \\ \cline{2-4}
\raisebox{1.5ex}[0pt]{712} & &  & \\ \hline
\multicolumn{4}{c}{$S_0$}
\end{tabular} $\Rightarrow$
\begin{tabular}{|c|c|c|c|}
\hline
 & 45 & 56 & 67 \\ \cline{2-4}
\raisebox{1.5ex}[0pt]{123} & $\bullet$ & $\bullet$ & $\times$ \\ \hline
 & 56 & 67 & 71 \\ \cline{2-4}
\raisebox{1.5ex}[0pt]{234} & &  &  \\ \hline
 & 67 & 71 & 12 \\ \cline{2-4}
\raisebox{1.5ex}[0pt]{345} & &  & $\times$ \\ \hline
 & 71 & 12 & 23 \\ \cline{2-4}
\raisebox{1.5ex}[0pt]{456} & $\bullet$ & $\times$ & $\times$ \\ \hline
 & 12 & 23 & 34 \\ \cline{2-4}
\raisebox{1.5ex}[0pt]{567} & $\times$ & $\bullet$ & $\times$\\ \hline
 & 23 & 34 & 45 \\ \cline{2-4}
\raisebox{1.5ex}[0pt]{671} & $\bullet$ &  &$\times$  \\ \hline
 & 34 & 45 & 56 \\ \cline{2-4}
\raisebox{1.5ex}[0pt]{712} & &  & \\ \hline
\multicolumn{4}{c}{$S^{\prime}_0$}
\end{tabular} $\Rightarrow$
\begin{tabular}{|c|c|c|c|}
\hline
 & 45 & 56 & 67 \\ \cline{2-4}
\raisebox{1.5ex}[0pt]{123} & $\bullet$ & $\bullet$ & $\times$ \\ \hline
 & 56 & 67 & 71 \\ \cline{2-4}
\raisebox{1.5ex}[0pt]{234} &  &  &   \\ \hline
 & 67 & 71 & 12 \\ \cline{2-4}
\raisebox{1.5ex}[0pt]{345} & &  & $\times$ \\ \hline
 & 71 & 12 & 23 \\ \cline{2-4}
\raisebox{1.5ex}[0pt]{456} & $\bullet$ & $\times$ & $\times$ \\ \hline
 & 12 & 23 & 34 \\ \cline{2-4}
\raisebox{1.5ex}[0pt]{567} & $\times$ & $\bullet$ &$\times$ \\ \hline
 & 23 & 34 & 45 \\ \cline{2-4}
\raisebox{1.5ex}[0pt]{671} & $\bullet$ &  & $\times$ \\ \hline
 & 34 & 45 & 56 \\ \cline{2-4}
\raisebox{1.5ex}[0pt]{712} & $\times$  & $\bullet$ & $\times$  \\ \hline
\multicolumn{4}{c}{$S^{\prime\prime}_0$}
\end{tabular}

\vspace{0.1cm}
\begin{tabular}{|c|c|c|c|}
\hline
 & 45 & 56 & 67 \\ \cline{2-4}
\raisebox{1.5ex}[0pt]{123} & $\bullet$ & $\bullet$ & $\times$ \\ \hline
 & 56 & 67 & 71 \\ \cline{2-4}
\raisebox{1.5ex}[0pt]{234} & $\bullet$ & $\times$ & $\times$  \\ \hline
 & 67 & 71 & 12 \\ \cline{2-4}
\raisebox{1.5ex}[0pt]{345} & &  & $\times$ \\ \hline
 & 71 & 12 & 23 \\ \cline{2-4}
\raisebox{1.5ex}[0pt]{456} & $\bullet$ & $\times$ & $\times$ \\ \hline
 & 12 & 23 & 34 \\ \cline{2-4}
\raisebox{1.5ex}[0pt]{567} & $\times$ & $\bullet$ &$\times$ \\ \hline
 & 23 & 34 & 45 \\ \cline{2-4}
\raisebox{1.5ex}[0pt]{671} & $\bullet$ &  & $\times$ \\ \hline
 & 34 & 45 & 56 \\ \cline{2-4}
\raisebox{1.5ex}[0pt]{712} & $\times$  & $\bullet$ & $\times$  \\ \hline
\multicolumn{4}{c}{$S_1$}
\end{tabular} $\Rightarrow$
\begin{tabular}{|c|c|c|c|}
\hline
 & 45 & 56 & 67 \\ \cline{2-4}
\raisebox{1.5ex}[0pt]{123} & $\bullet$ & $\bullet$ & $\times$ \\ \hline
 & 56 & 67 & 71 \\ \cline{2-4}
\raisebox{1.5ex}[0pt]{234} & $\bullet$ & $\times$ & $\times$  \\ \hline
 & 67 & 71 & 12 \\ \cline{2-4}
\raisebox{1.5ex}[0pt]{345} & $\times$ & $\bullet$ & $\times$ \\ \hline
 & 71 & 12 & 23 \\ \cline{2-4}
\raisebox{1.5ex}[0pt]{456} & $\bullet$ & $\times$ & $\times$ \\ \hline
 & 12 & 23 & 34 \\ \cline{2-4}
\raisebox{1.5ex}[0pt]{567} & $\times$ & $\bullet$ &$\times$ \\ \hline
 & 23 & 34 & 45 \\ \cline{2-4}
\raisebox{1.5ex}[0pt]{671} & $\bullet$ &  & $\times$ \\ \hline
 & 34 & 45 & 56 \\ \cline{2-4}
\raisebox{1.5ex}[0pt]{712} & $\times$  & $\bullet$ & $\times$  \\ \hline
\multicolumn{4}{c}{$S^{\prime}_1$}
\end{tabular}\hspace{0.1cm}
\begin{tabular}{|c|c|c|c|}
\hline
 & 45 & 56 & 67 \\ \cline{2-4}
\raisebox{1.5ex}[0pt]{123} & $\bullet$ & $\bullet$ & $\times$ \\ \hline
 & 56 & 67 & 71 \\ \cline{2-4}
\raisebox{1.5ex}[0pt]{234} & $\bullet$ & $\times$ & $\times$  \\ \hline
 & 67 & 71 & 12 \\ \cline{2-4}
\raisebox{1.5ex}[0pt]{345} & $\times$ & $\bullet$ & $\times$ \\ \hline
 & 71 & 12 & 23 \\ \cline{2-4}
\raisebox{1.5ex}[0pt]{456} & $\bullet$ & $\times$ & $\times$ \\ \hline
 & 12 & 23 & 34 \\ \cline{2-4}
\raisebox{1.5ex}[0pt]{567} & $\times$ & $\bullet$ &$\times$ \\ \hline
 & 23 & 34 & 45 \\ \cline{2-4}
\raisebox{1.5ex}[0pt]{671} & $\bullet$ & $\times$  & $\times$ \\ \hline
 & 34 & 45 & 56 \\ \cline{2-4}
\raisebox{1.5ex}[0pt]{712} & $\times$  & $\bullet$ & $\times$  \\ \hline
\multicolumn{4}{c}{$S^{\prime}_{1-1}$}
\end{tabular} \hspace{0.1cm}
\begin{tabular}{|c|c|c|c|}
\hline
 & 45 & 56 & 67 \\ \cline{2-4}
\raisebox{1.5ex}[0pt]{123} & $\bullet$ & $\bullet$ & $\times$ \\ \hline
 & 56 & 67 & 71 \\ \cline{2-4}
\raisebox{1.5ex}[0pt]{234} & $\bullet$ & $\times$ & $\times$  \\ \hline
 & 67 & 71 & 12 \\ \cline{2-4}
\raisebox{1.5ex}[0pt]{345} & $\times$ & $\bullet$ & $\times$ \\ \hline
 & 71 & 12 & 23 \\ \cline{2-4}
\raisebox{1.5ex}[0pt]{456} & $\bullet$ & $\times$ & $\times$ \\ \hline
 & 12 & 23 & 34 \\ \cline{2-4}
\raisebox{1.5ex}[0pt]{567} & $\times$ & $\bullet$ &$\times$ \\ \hline
 & 23 & 34 & 45 \\ \cline{2-4}
\raisebox{1.5ex}[0pt]{671} & $\bullet$ & $\bullet$  & $\times$ \\ \hline
 & 34 & 45 & 56 \\ \cline{2-4}
\raisebox{1.5ex}[0pt]{712} & $\times$  & $\bullet$ & $\times$  \\ \hline
\multicolumn{4}{c}{$S^{\prime}_{1-2}$}
\end{tabular}

\vspace{0.1cm}
\begin{tabular}{|c|c|c|c|}
\hline
 & 45 & 56 & 67 \\ \cline{2-4}
\raisebox{1.5ex}[0pt]{123} & $\bullet$ & $\bullet$ & $\times$ \\ \hline
 & 56 & 67 & 71 \\ \cline{2-4}
\raisebox{1.5ex}[0pt]{234} & $\times$ & $\bullet$ & $\times$  \\ \hline
 & 67 & 71 & 12 \\ \cline{2-4}
\raisebox{1.5ex}[0pt]{345} & &  & $\times$ \\ \hline
 & 71 & 12 & 23 \\ \cline{2-4}
\raisebox{1.5ex}[0pt]{456} & $\bullet$ & $\times$ & $\times$ \\ \hline
 & 12 & 23 & 34 \\ \cline{2-4}
\raisebox{1.5ex}[0pt]{567} & $\times$ & $\bullet$ &$\times$ \\ \hline
 & 23 & 34 & 45 \\ \cline{2-4}
\raisebox{1.5ex}[0pt]{671} & $\bullet$ &  & $\times$ \\ \hline
 & 34 & 45 & 56 \\ \cline{2-4}
\raisebox{1.5ex}[0pt]{712} & $\times$  & $\bullet$ & $\times$  \\ \hline
\multicolumn{4}{c}{$S_2$}
\end{tabular} $\Rightarrow$
\begin{tabular}{|c|c|c|c|}
\hline
 & 45 & 56 & 67 \\ \cline{2-4}
\raisebox{1.5ex}[0pt]{123} & $\bullet$ & $\bullet$ & $\times$ \\ \hline
 & 56 & 67 & 71 \\ \cline{2-4}
\raisebox{1.5ex}[0pt]{234} & $\times$ & $\bullet$ & $\times$  \\ \hline
 & 67 & 71 & 12 \\ \cline{2-4}
\raisebox{1.5ex}[0pt]{345} & &  & $\times$ \\ \hline
 & 71 & 12 & 23 \\ \cline{2-4}
\raisebox{1.5ex}[0pt]{456} & $\bullet$ & $\times$ & $\times$ \\ \hline
 & 12 & 23 & 34 \\ \cline{2-4}
\raisebox{1.5ex}[0pt]{567} & $\times$ & $\bullet$ &$\times$ \\ \hline
 & 23 & 34 & 45 \\ \cline{2-4}
\raisebox{1.5ex}[0pt]{671} & $\bullet$ & $\times$ & $\times$ \\ \hline
 & 34 & 45 & 56 \\ \cline{2-4}
\raisebox{1.5ex}[0pt]{712} & $\times$  & $\bullet$ & $\times$  \\ \hline
\multicolumn{4}{c}{$S_2^{\prime}$}
\end{tabular} \hspace{0.1cm}
\begin{tabular}{|c|c|c|c|}
\hline
 & 45 & 56 & 67 \\ \cline{2-4}
\raisebox{1.5ex}[0pt]{123} & $\bullet$ & $\bullet$ & $\times$ \\ \hline
 & 56 & 67 & 71 \\ \cline{2-4}
\raisebox{1.5ex}[0pt]{234} & $\times$ & $\bullet$ & $\times$  \\ \hline
 & 67 & 71 & 12 \\ \cline{2-4}
\raisebox{1.5ex}[0pt]{345} &$\bullet$ & $\bullet$ & $\times$ \\ \hline
 & 71 & 12 & 23 \\ \cline{2-4}
\raisebox{1.5ex}[0pt]{456} & $\bullet$ & $\times$ & $\times$ \\ \hline
 & 12 & 23 & 34 \\ \cline{2-4}
\raisebox{1.5ex}[0pt]{567} & $\times$ & $\bullet$ &$\times$ \\ \hline
 & 23 & 34 & 45 \\ \cline{2-4}
\raisebox{1.5ex}[0pt]{671} & $\bullet$ & $\times$ & $\times$ \\ \hline
 & 34 & 45 & 56 \\ \cline{2-4}
\raisebox{1.5ex}[0pt]{712} & $\times$  & $\bullet$ & $\times$  \\ \hline
\multicolumn{4}{c}{$S_{2-1}^{\prime}$}
\end{tabular} \hspace{0.1cm}\begin{tabular}{|c|c|c|c|}
\hline
 & 45 & 56 & 67 \\ \cline{2-4}
\raisebox{1.5ex}[0pt]{123} & $\bullet$ & $\bullet$ & $\times$ \\ \hline
 & 56 & 67 & 71 \\ \cline{2-4}
\raisebox{1.5ex}[0pt]{234} & $\times$ & $\bullet$ & $\times$  \\ \hline
 & 67 & 71 & 12 \\ \cline{2-4}
\raisebox{1.5ex}[0pt]{345} & $\times$& $\bullet$ & $\times$ \\ \hline
 & 71 & 12 & 23 \\ \cline{2-4}
\raisebox{1.5ex}[0pt]{456} & $\bullet$ & $\times$ & $\times$ \\ \hline
 & 12 & 23 & 34 \\ \cline{2-4}
\raisebox{1.5ex}[0pt]{567} & $\times$ & $\bullet$ &$\times$ \\ \hline
 & 23 & 34 & 45 \\ \cline{2-4}
\raisebox{1.5ex}[0pt]{671} & $\bullet$ & $\times$ & $\times$ \\ \hline
 & 34 & 45 & 56 \\ \cline{2-4}
\raisebox{1.5ex}[0pt]{712} & $\times$  & $\bullet$ & $\times$  \\ \hline
\multicolumn{4}{c}{$S_{2-2}^{\prime}$}
\end{tabular} \hspace{0.1cm}
}
\caption{ }\label{table-1}
\end{table}

\begin{figure}
\centerline{\epsfbox{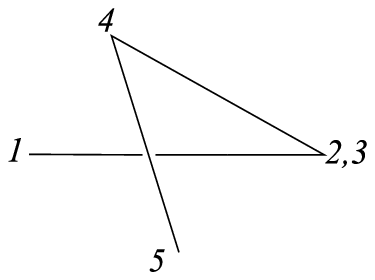}}
\caption{ }
\label{fig8}
\end{figure}

\vspace{0.1cm}
Now we prove the ``if'' part of the theorem. Suppose $P$ is a heptagonal knot satisfying {\em RS-I}, {\em II} or {\em III}. Let $N$ be a plane orthogonal to $\overrightarrow{23}$, and $\pi:\mathbb{R}^3 \equiv N\times\mathbb{R} \rightarrow N$ be the orthogonal projection onto $N$ such that the vertex $3$ is above the vertex $2$ with respect to the $\mathbb{R}$-coordinate. We will construct a diagram of $P$ from the projected image $\pi(P)$. Without loss of generality it can be assumed that $\epsilon(123,45)=1$ and $\epsilon(123,56)=-1$. Then, since the vertex $3$ is above $2$, the edge $e_{45}$ should pass above $e_{12}$ as illustrated in Figure \ref{fig8}.

\begin{figure}
\centerline{\epsfbox{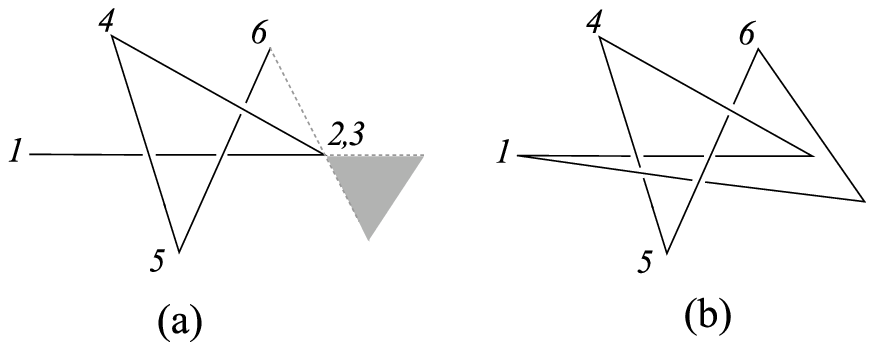}}
\caption{ }
\label{fig9}
\end{figure}

Suppose $P$ corresponds to {\em RS-I}. Then similarly $e_{56}$ passes above $e_{12}$ and below $e_{34}$. Note that if $7$ belongs to $H_{123}^-$, then $e_{23}$ can not penetrate $\Delta_{671}$. Hence $7 \in H_{123}^+$. Since $\epsilon(567,23)=-1$, the point $\pi(e_{23})$ should belong to $\pi(\Delta_{567})$. Therefore the point $\pi(7)$ belongs to the shaded region shown in Figure \ref{fig9}-(a). From this we know that $\pi(7) \notin \pi(\Delta_{456})$. And clearly $\pi(1) \notin \pi(\Delta_{456})$. Hence, for $\epsilon(456,71)$ to be nonzero, $\pi(e_{71})$ should intersect both $\pi(e_{45})$ and $\pi(e_{56})$. In fact $e_{71}$ passes above $e_{45}$ because $\epsilon(712,45)=1$ and $e_{45}$ passes above $e_{12}$. Therefore, since $\epsilon(456,71)$ is nonzero, $e_{71}$ should pass below $e_{56}$. The resulting diagram represents figure-8 as seen in Figure \ref{fig9}-(b). Also when $P$ corresponds to {\em RS-II}, we can obtain a diagram of figure-8 in the same way.

\begin{figure}
\centerline{\epsfbox{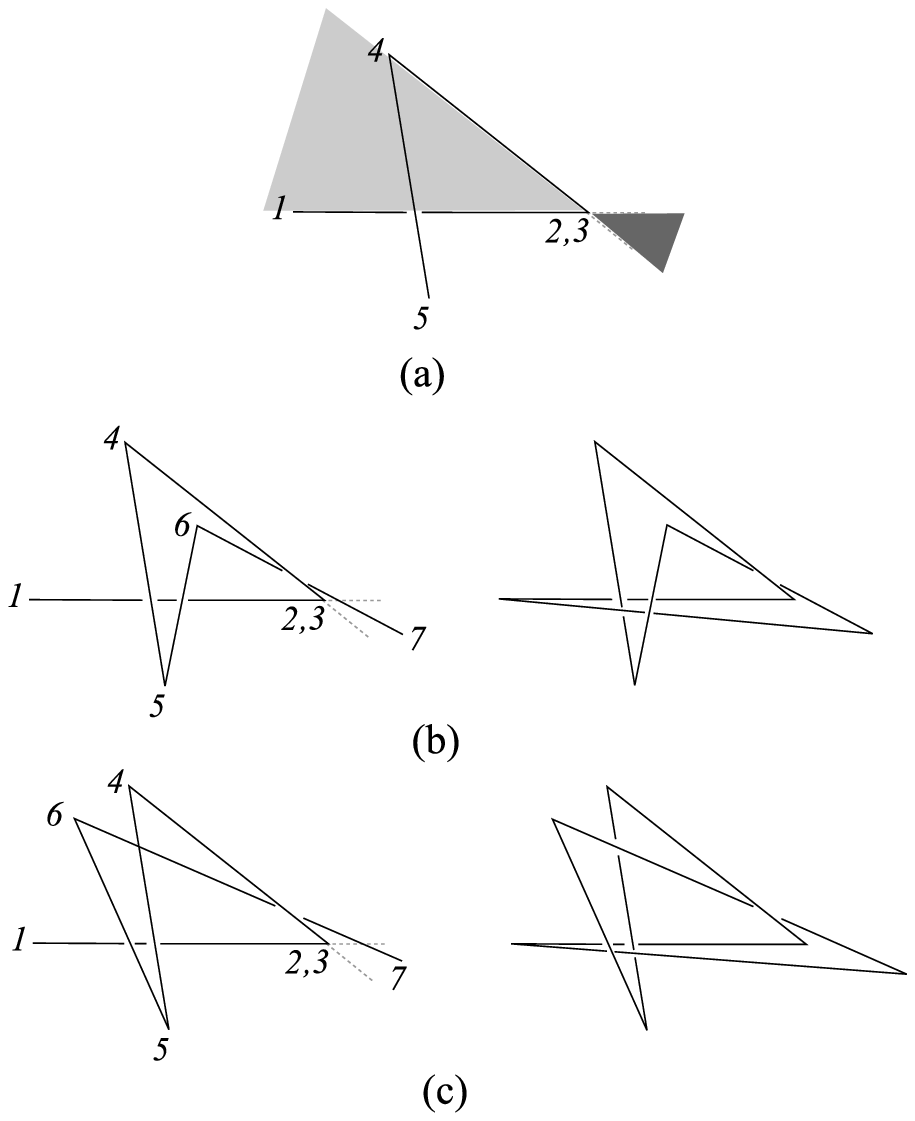}}
\caption{ }
\label{fig10}
\end{figure}

Suppose $P$ corresponds to {\em RS-III}. Again assume  $\epsilon(123,45)=1$ and $\epsilon(123,56) =-1$.
Then also in this case we have that $6\in H_{123}^-$ and $7 \in H_{123}^+$. Especially it should hold that $6 \in H_{234}^+$ and $7\in H_{234}^-$, because $\epsilon(234,67)=-1$. Hence the vertex $6$ is projected into the shaded region in the top-left of Figure \ref{fig10}-(a) and the vertex $7$ into the bottom-right. Now we observe two possible cases according to the position of $\pi(6)$ with respect to $\pi(e_{45})$ as shown in Figure \ref{fig10}-(b) and (c). Again since $\epsilon(234,67)$ is nonzero, $e_{67}$ should pass below $e_{34}$ in both diagrams. As discussed in the case of {\em RS-I}, from $\epsilon(456,71)=1$, we know that $e_{71}$ passes above $e_{45}$ and below $e_{56}$. Then the resulting diagram in (b) represents figure-8. In (c), for $e_{71}$ to pass above $e_{45}$ and below $e_{56}$, the two vertices $6$ and $7$ should belong to $H_{145}^-$, which implies that $e_{67}$ passes above $e_{45}$. Therefore also the resulting diagram in (c) represents figure-8.


\end{document}